\documentclass[twocolumn,preprintnumbers,amsmath,amssymb]{revtex4}

\def\aslash{\raisebox{0.40ex}{\tiny$\setminus$}}
\def\slash{\raisebox{0.40ex}{\tiny/}}

\begin{document}

\title{GENERALIZATION OF NUMERICAL SERIES AND ITS RELATIONSHIP WITH THE POLYNOMIAL EQUATIONS AND ARITHMETIC TRAPEZOIDS}

\author{Victor Enrique Vizcarra Ruiz}

\affiliation{Departamento de F\'isica\\Universidade Estadual de Maring\'a\\
Av. Colombo, 5790, Maring\'a - Pr, 87020-900, Brasil.}

\date{\today}

\begin{abstract}
{\indent The close relationship among the polynomial functions and Fibonacci numerical sequences is shown in this paper. These numerical sequences are defined by the recurrence equation $x_{k + n} = \displaystyle\sum_{j = 0}^{n-1}\alpha_j x_{k + j}$, where $n$ is the polynomial degree and $\alpha$'s, the polynomial coefficients. The arithmetic trapezoid resulting from the recurrence equations is also shown. This trapezoid is nothing but a generalization of Pascal's Triangle. Trapezoid is a convenient name because the form it appears does not have the `upper end` of a usual triangle. This study shows that each polynomial generates infinite sequences, and that each sequence generates only a single arithmetic trapezoid.}
\end{abstract}

\maketitle

\section{Introduction}

Although known by Indian poets over the centuries before Christ \cite{pingala,singh}, the sequence $1,1,2,3,5,8,13, \dots$, was only introduced in the West in 1202 by Leonardo Bigollo, referred to as Fibonacci \cite{fibonacci}, in his famous problem of rabbit breeding. Since then, many researchers have dedicated their studies to explore the properties of this sequence, and even to create new series, as Lucas \cite{lucas} and Pell\footnote{ In spite of the fact that the sequence is named after Pell, it seems that the true creator was William Brouncker.}, with new properties and relating them to the Fibonacci sequence. There is a wide effort to relate all these series with some general property. This study shows how to relate them all, in addition to the possibility of creating new series, infinite series, from a sequence generalization. The Fibonacci sequence, that of Lucas, Pell's one, and many others, are just special cases of this generalization, and the starting point will be its relationship with the golden number $\phi = 1.6180339\dots$. This number is the value to which the ratio of two consecutive numbers of the Fibonacci sequence converges, that is, $x_{k+1}/x_k \approx \phi$, where $x_k$ is any number of the sequence, and $k$ is the index which indicates the location of that number in the sequence. The higher the value of $k$, the more the ratio of these two numbers is closer to the golden value. Since this golden number is the solution of a particular quadratic equation, it is possible that other particular solutions are also values to which the ratio of two consecutive numbers of other sequences converges. Therefore, it means that the general quadratic equation is related to a general recurrence equation, which is able to reproduce the numbers of the Fibonacci series. The series itself is defined by the recurrence relation $x_{k+2} = x_{k+1} + x_k$, which was also known by the Indians. Binet \cite{binet} obtained a formula that directly calculated any number of the sequence, without using the recurrence equation, that is, $x_k = (\phi^n - (1 - \phi)^n)/\sqrt{5}$. These two equations will be generalized in this study, as well as their corresponding golden numbers from the polynomial equations.

Another important generalization is Pascal's arithmetical triangle because it is related to the Fibonacci sequence; in the sense that the sum of certain diagonals of the triangle reproduces the sequence. This study shows two important points. First, that the arithmetical triangle is not a triangle, but a trapezoid, since the form of its construction does not have the `upper end' of a usual triangle. Pascal's triangle, known by ancient people long before Pascal, is just a particular case of this \textit{arithmetical trapezoid}. And that is how it will be called from now on. And second, that for each numerical sequence there is a single arithmetic trapezoid. After this study Pascal's triangle will most likely be called Pascal's trapezoid.

How the number sequences are related to both the polynomial equations and the arithmetical trapezoid will be formally shown in the following topic.

\section{Generalization}

\subsection{Relationship with the polynomial equations}

The numerical sequences are related to the roots of the polynomial equations through the coefficients of such polynomials. In the case of the \textit{\textbf{quadratic equation}} $ax^2 + bx + c = 0$, where $a, b, c \in \mathbb{R}$, the general recurrence relationship associated with this polynomial is given as $ax_{k+2} + bx_{k+1} + cx_k = 0$. However, without loss of generality, it can be given as it follows,
\begin{equation}
\left\{
\begin{array}{l}
x^2 = \alpha x + \beta,\\[0.1cm]
x_{k+2} = \alpha x_{k+1} + \beta x_k
\end{array}
\right.\label{segundograu}
\end{equation}
where $\alpha = -b/a$ ans $\beta = -c/a$, with $\alpha, \beta \in \mathbb{R}$. The initial values $x_0$ and $x_1$ are the recurrence relation seeds and can admit complex values. The index $k$ is, again, the position of the number $x_k$ in the sequence. The coefficients $\alpha$ and $\beta$ of the polynomial equation in (\ref{segundograu}), whose roots are clearly $\phi = (\alpha + \sigma)/2$ and $\varphi = (\alpha - \sigma)/2$, where $\sigma = \sqrt{\alpha^2 + 4\beta}$, are the parameters that generalize the number sequence $x_0, x_1, x_2, x_3, x_4, x_5, \cdots$, that is,
\begin{widetext}
\begin{equation}
\begin{array}{ccccc}
\underbrace{x_0},& \underbrace{x_1},& \underbrace{\alpha x_1 + \beta x_0},& \underbrace{(\alpha^2 + \beta) x_1 + \alpha\beta x_0},&\hspace{-.5cm} \underbrace{\alpha(\alpha^2 + 2\beta)x_1 + \beta(\alpha^2 + \beta)x_0},\\[.3cm]
x_0&x_1&x_2&x_3&\hspace{-.5cm}x_4 \\[.3cm]
&&&& \hspace{-1.5cm}\underbrace{(\alpha^4 + 3\alpha^2\beta + \beta^2)x_1 + (\alpha^3\beta + 2\alpha\beta^2)x_0},\ \cdots\\
&&&&\hspace{-1.5cm} x_5\ \ \ \ \
\end{array}
\end{equation}
\end{widetext}
as well as the equation that directly calculates any number of the sequence,
\begin{equation}
x_k = \frac{(\phi - \alpha)x_0 + x_1}{\sqrt{\alpha^2 + 4\beta}}(\phi^k - \varphi^k) + \varphi^kx_0.
\label{xn2grau}
\end{equation}

The above equation  (\ref{xn2grau}), s Binet?s general formula. For the particular cases in which the coefficients $\alpha$ and $\beta$ correspond to unit values, and the seeds are equal to $(x_0, x_1) = (0, 1)$ and the $(x_0, x_1) = (2, 1)$, there will exactly be the Fibonacci and Lucas? sequences, respectively; that is,
\begin{widetext}
\begin{equation}
\begin{array}{c}
\alpha = 1,\ \beta = 1 ;\ \ \ \ \
\left\{
\begin{array}{c}
\displaystyle x^2 - x - 1 = 0 \\
\displaystyle x_{k+2} = x_{k+1} + x_k \\
\end{array}\right.;\ \ \ \ \ \phi = \displaystyle\frac{1 + \sqrt{5}}2,\ \  \varphi = \displaystyle\frac{1 - \sqrt{5}}2
\\[0.5cm]
\begin{array}{cc}
\textrm{Fibonacci} \ \ \ \ \ \ &\ \ \ \ \textrm{Lucas} \\[0.2cm]
x_0 = 0,\ \ x_1 = 1\ \ \ \ \ \ \ &\ \ \ \ x_0 = 2,\ \ x_1 = 1 \\[0.1cm]
0, 1, 1, 2, 3, 5, \cdots\ \ \ \ \ \ &\ \ \ \ \ \ 2, 1, 3, 4, 7, 11, \dots \\[0.1cm]
x_k = \displaystyle\frac{\phi^k - \varphi^k}{\sqrt{5}}\ \ \ \ &\ \ \ \ x_k = \phi^k + \varphi^k \\
\end{array}
\end{array}\nonumber
\end{equation}
\end{widetext}

The natural consequence of this generalization is that infinite number sequences can be obtained, varying not only the seeds values $(x_0, x_1)$, but also the values of the coefficients $\alpha$ and $\beta$. As $\alpha$ and $\beta$ are real values, in some cases the roots $\phi$ and $\varphi$ will be complex quantities.

The root $\phi$ is the convergence value of fractions between two successive consecutive numbers of the sequence, that is,
$$\lim_{k \rightarrow \infty} \frac{x_{k+1}}{x_k} = \phi,$$ e que pode ser verificado utilizando a rela\c c\~ao (\ref{segundograu}), isto \',
\begin{equation}
\left\{
\begin{array}{l}
\phi = \alpha + \displaystyle\frac{\beta}{\phi},\\[0.3cm]
\displaystyle\lim_{k \rightarrow \infty} \displaystyle\frac{x_{k+2}}{x_{k+1}} = \alpha +  \displaystyle\frac{\beta}{\displaystyle\lim_{k \rightarrow \infty} x_{k+1}/x_k}
\end{array}
\right.
\end{equation}
IIn the literature, the symbol $\phi$ is specifically the rational number 1.618033989$\dots$, which is called golden number. The same symbol will be used to represent the root $(\alpha + \sigma)/2$ of the polynomial, as it is depicted in Eq. (\ref{raizes2grau}), and it will be  the \textit{`sequence golden number'}. The golden number for this particular case $(\alpha, \beta) = (2, 1)$, for example, corresponds to the already known silver number $\phi = 1 + \sqrt{2}$, which is the characteristic of all Pell's sequences and that depends on the seed values $x_0$ e $x_1$. Aiming practicability, in some cases the symbol $\phi(\alpha,\beta)$ is used in order to refer to the same golden number of Eq. (\ref{raizes2grau}). Considering the silver number, for example, it will be written $\phi(2,1) = 2.414213562\dots$ The particular set in which the polynomial coefficients are specifically unit values, that is, $(\alpha, \beta) = (1, 1)$, and only referring to this set, will be called `Fibonacci' or `golden'; and the set in which $(\alpha, \beta) = (2, 1)$, will be called `silver'. In the literature there are already names of the metallic families \cite{metalicas}, for the different values of $\alpha$ and $\beta$. The purpose of such names is just mere classification. The term \textit{set} is suitable because it gathers a family dependent on the seed values $(x_0, x_1)$. Fibonacci and Lucas' sequences are part of the Fibonacci sets. Pell's sequence corresponds to the silver group, whose seeds are $(x_0,x_1) = (0,1)$.

In the case of the \textit{\textbf{cubic equation}}, $ax^3 + bx^2 + cx + d = 0$, where $a, b, c \in \mathbb{R}$, and $\alpha, \beta, \gamma \in \mathbb{R}$, with $\alpha = -b/a$, $\beta = -c/a$ and $\gamma = -d/a$, the polynomial and recurrence equation associated to it is given as,
\begin{equation}
\left\{
\begin{array}{l}
x^3 = \alpha x^2 + \beta x + \gamma,\\[0.1cm]
x_{k+3} = \alpha x_{k+2} + \beta x_{k+1} + \gamma x_k,
\end{array}
\right.\label{terceirograu}
\end{equation}
where $x_0$, $x_1$ and $x_2$ are the seed values. The golden numbers, or roots, associated to the sequences from equation (\ref{terceirograu}), are also given as polynomial coefficients, i.e., $\phi = (\alpha + \sigma_1 + \sigma_2)/3$, $\varphi = (\alpha\ \aslash\ \sigma_1\ \slash\ \sigma_2)/3$ and $\psi = (\alpha\ \slash\ \sigma_1\ \aslash\ \sigma_2)/3$, where $\sigma_1 = \left[\left(A + \sqrt{A^2 - 4B^3}\right)/2\right]^{1/3}$, $\sigma_2 = \left[\left(A - \sqrt{A^2 - 4B^3}\right)/2\right]^{1/3}$ and where $A = 2\alpha^3 + 9\alpha\beta + 27\gamma$ and $B = \alpha^2 + 3\beta$.

Two new typographical marks that represent arithmetical symbols are introduced, `\slash' and `\,\aslash', slash and backlash, respectively, in order to simplify the mathematical writing, making it clear that they are only complex numbers that `imitate' operators, such as  `+' and `--'. Thus, $\slash1 = (-1 + i\sqrt{3})/2$ and $\aslash1 = (-1 - i\sqrt{3})/2$ is written, with $i = \sqrt{-1}$. Therefore, the relationships between golden numbers there will be $\phi + \varphi + \psi = \alpha$, $\phi\ \slash\ \varphi\ \aslash\ \psi = \sigma_1$ and $\phi\ \aslash\ \varphi\ \slash\ \psi = \sigma_2$.

With the introduction of `\,$\slash$\,' and `\,$\aslash$\,' new mathematical structures with intrinsic properties might emerge, such as $\phi\ \slash\ \phi\ \aslash\ \phi = 0$, which can be applied to the strange properties of the universe of quarks. However, this will be subject for a further study.

The equation which directly calculates an arbitrary number of the sequence (\ref{terceirograu}), Binet's generalized equation, is given by,

\begin{widetext}
\begin{equation}
\begin{array}{l}
\hspace{-1.5cm}x_k = \displaystyle\left[\frac{(\alpha -\varphi - \psi - 1)x_2 + (\beta + \varphi\psi + \varphi + \psi)x_1 + (\gamma - \varphi\psi) x_0}{(\phi - \varphi)(\varphi - \psi)(\phi - \psi)}\right]\left(\frac{\psi - \varphi}{\phi - 1}\phi^k - \frac{\psi - \phi}{\varphi - 1}\varphi^k + \frac{\varphi - \phi}{\psi - 1}\psi^k\right)\\[0.5cm]
\hspace{6.5cm}+ \displaystyle\frac{x_2 - (\varphi + 1)x_1 + \varphi x_0}{(\psi - 1)(\psi - \varphi)}\psi^k
- \displaystyle\frac{x_2 - (\psi + 1)x_1 + \psi x_0}{(\varphi - 1)(\psi - \varphi)}\varphi^k
\end{array}
\label{xn3grau}
\end{equation}
\end{widetext}
The numbers obtained in the above equation (\ref{xn3grau}), or by its recurrence equation, given in (\ref{terceirograu}), have the basic property, as well as in the quadratic case: the ratio of two consecutive numbers converge to their golden number, that is, $\phi = \lim_{k \rightarrow \infty} x_{k + 1}/x_k$. This property can be verified in the relationships (\ref{terceirograu}),
\begin{widetext}
\begin{equation}
\left\{
\begin{array}{l}
\phi = \alpha + \displaystyle\frac{\beta}{\phi} + \displaystyle\frac{\gamma}{\phi^2},\\[0.3cm]
\displaystyle\lim_{k \rightarrow \infty}\frac{x_{k+3}}{x_{k+2}} = \alpha + \frac{\beta}{\displaystyle\lim_{k \rightarrow \infty} x_{k+2}/x_{k+1}} + \frac{\gamma}{\left(\displaystyle\lim_{k \rightarrow \infty} x_{k+2}/x_{k+1}\right)
\left(\displaystyle\lim_{k \rightarrow \infty} x_{k+1}/x_k\right)}
\end{array}
\right.
\end{equation}
\end{widetext}
However, it is not so simple to verify this convergence for the numbers $\varphi$ and $\psi$ because these two numbers are not pure real values, but complex ones. This can be seen in the expressions for $\varphi$ and $\psi$, exactly on the pseudo-signs slash `\,\slash\,' and backslash `\,\aslash\,'. Considering the fact that these pseudo-signals represent the complex numbers $(-1 + i\sqrt{3})/2$ and $(-1 - i\sqrt{3})/2$, respectively, the roots $\varphi$ and $\psi$ can be written in function of $\phi$, and, thus, to verify the convergence of the numbers of the series for $\varphi$ and $\psi$, that is,
\begin{widetext}
\begin{equation}\left\{
\begin{array}{ccc}
\varphi & = & \displaystyle\frac12\left[
\alpha - \lim_{k \rightarrow \infty} \frac{x_{k + 1}}{x_k}
 + \sqrt{\left(\alpha - \lim_{k \rightarrow \infty} \frac{x_{k + 1}}{x_k}\right)^2 - 4\gamma\lim_{k \rightarrow \infty} \frac{x_k}{x_{k + 1}}}\right]\nonumber\\[0.8cm]
\psi & = & \displaystyle\frac12\left[\alpha - \lim_{k \rightarrow \infty} \frac{x_{k + 1}}{x_k}
 - \sqrt{\left(\alpha - \lim_{k \rightarrow \infty} \frac{x_{k + 1}}{x_k}\right)^2 - 4\gamma\lim_{k \rightarrow \infty} \frac{x_k}{x_{k + 1}}}\right]\nonumber
\end{array}\right.
\end{equation}
\end{widetext}
Numerically, the roots for this Fibonacci set $(\alpha,\beta,\gamma) = (1,1,1)$, with the seeds $(x_0,x_1,x_2) = (0,1,1)$, are written as $\phi = 1,839286755\dots$, sometime know as the tribonacci constant, $\varphi = -0,4196433776\dots$$ + i\,0,6062907292\dots$ and $\psi = -0,4196433776\dots$$ - i\,0,6062907292\dots$. These are the values to which the ratios of two successive numbers of the series converge.\\

When generalizing the basic equations of a \textit{\textbf{polynomial of any $\boldsymbol{n}$-degree}}, it is given as,
\begin{equation}
\left\{
\begin{array}{l}
x^n = \displaystyle\sum_{j = 0}^{n-1}\alpha_j x^j,\\[0.5cm]
x_{k + n} = \displaystyle\sum_{j = 0}^{n-1}\alpha_j x_{k + j}.
\end{array}
\right.\label{graugeral}
\end{equation}
where the first expression corresponds to the polynomial equation of $n$th degree, and the second, the recurrence equation of the numeric sequence, resulting from the polynomial equation, where  $\alpha_j\textrm{'s} \in \mathbb{R}$, are the polynomial coefficients. The roots of the general polynomial in (\ref{graugeral}), can be obtained by using any computer program that works with algebraic analysis. On the other hand, Binet's generalized equation, that is, the equation that calculates an arbitrary number of the sequence associated with this polynomial, is obtained due to the fact that the recurrence equation is linear and homogeneous. In this case, the chosen number is written as a linear combination of all the solutions of the polynomial equation, that is, of their roots, as it follows,
\begin{equation}
x_k = \sum_{j = 1}^nw_j\phi_j^k + w_{n+1}.
\label{equacaolinear}
\end{equation}
where $\phi$'s represents the roots, and $w$'s, the coefficients that will be determined. The $n$ index represents the degree of the polynomial. Since there are $n+1$ unknowns $w$'s, linear equations will be written, adjusting them to the first numbers of the sequence $n+1$, that is,
\begin{equation}\left[
\begin{array}{cccccc}
1 & 1 & 1 & \dots & 1 & 1 \\[0.1cm]
\phi_1 & \phi_2 & \phi_3 & \dots & \phi_n & 1 \\[0.1cm]
\phi_1^2 & \phi_2^2 & \phi_3^2 & \dots & \phi_n^2 & 1 \\[0.1cm]
 \vdots &  &  & \vdots &  & \vdots \\[0.1cm]
\phi_1^n & \phi_2^n & \phi_3^n & \dots & \phi_n^n & 1 \\
\end{array}\right]
\left[
\begin{array}{c}
w_1 \\[0.1cm]
w_2 \\[0.1cm]
w_3 \\[0.1cm]
\vdots \\[0.1cm]
w_{n+1} \\
\end{array}\right] = 
\left[
\begin{array}{c}
x_0 \\[0.1cm]
x_1 \\[0.1cm]
x_2 \\[0.1cm]
\vdots \\[0.1cm]
x_n \\
\end{array}\right]\label{sistemageral}
\end{equation}
Among the results that might be obtained $w_{n+1} = 0$ will be found. The other findings will be related to each other. In this study, all the $w$'s are related to $w_1$, that is, $w_2 \propto w_1$, $w_3 \propto w_1$, $\hdots$, $w_n \propto w_1$. The reason is just to simplify the final expression. This is how the expressions (\ref{xn2grau}) and (\ref{xn3grau}) were obtained.

The expressions contained in (\ref{graugeral}), connect the numerical sequences with the polynomial equations. The fundamental property of the recurrence law that governs this numerical sequence is that given two consecutive numbers in this sequence, the ratio of these two numbers converges to their golden number $\phi$, that is, $\lim_{k \rightarrow \infty} (x_{k + n})/(x_{k + n - 1}) = \phi$.

Equation (\ref{equacaolinear}), whose coefficients $w$'s were determined by using the linear equation system in (\ref{sistemageral}), allows us to directly calculate any number of the sequence just by knowing its location $k$ within the sequence. This is Binet's equation principle.\\

Finally, a little about the roots of the polynomial equation will be seen. Given a polynomial of $n$-degree, its roots $\phi_1$, $\phi_2$, $\ldots$, $\phi_n$, are related to each other through the elementary symmetric polynomials, that is,
\begin{equation}
\left\{\begin{array}{ccc}
\phi_1 + \phi_2 + \phi_3 + \cdots + \phi_n & = & \ \ \alpha_{n - 1}\\[0.2cm]
\phi_1\phi_2 + \phi_1\phi_3 + \cdots + \phi_{n - 1}\phi_n & = & -\alpha_{n - 2}\\[0.2cm]
\cdots &  & \\[0.2cm]
\phi_1\phi_2\phi_3\cdots\phi_n & = & (-1)^{n - 1}\alpha_0
\end{array}\right.\label{polinomiossimetricos}
\end{equation}
As these roots are the golden numbers related to the numerical sequences that result from the polynomial equation, all the golden numbers of a particular group of numerical sequences associated with the polynomial coefficients $\alpha_0$, $\alpha_1$, $\ldots$, $\alpha_n$, are naturally related by these symmetric elementary polynomials (\ref{polinomiossimetricos}). Therefore, for example, when using the symbols of the expressions (\ref{segundograu}) and (\ref{terceirograu}), that is, $\phi = \phi_1$, $\varphi = \phi_2$, $\alpha = \alpha_1$ and $\beta = \alpha_0$, and, in case of the quadratic equation, $n = 2$, there will be $\phi + \varphi = \alpha$, $\phi\varphi = -\beta$, and $\phi - \varphi = \sigma$.

Similarly, for a cubic polynomial $\phi\, +\, \varphi + \psi = \alpha$, $\phi\varphi + \phi\psi + \varphi\psi = -\beta$ and $\phi\varphi\psi = \gamma$, and $\phi\ \slash\ \varphi\ \aslash\ \psi = \sigma_1$ and $\phi\ \aslash\ \varphi\ \slash\ \psi = \sigma_2$, are found.

Another point about these roots will be considered. Since these numbers are associated to the numerical sequences, through the convergence of the ratios of two successive consecutive numbers, it is necessary that all other roots, $\varphi$, $\psi$, \dots, are only written in terms of  $\phi = \lim_{k \rightarrow \infty} x_{k + 1}/x_k$. Only this way the convergence for these values will be verified.

Finally, the properties of each golden number are obtained by using the polynomial own equations. For example, for all Fibonacci sets of the quadratic case, the square of their golden number is the golden number itself by adding one unit, that is, $\phi^2 = \phi + 1$; and its inverse is itself by subtracting one unit, that is, $1/\phi = \phi - 1$. These properties, and many others, are generated from its general equation, $x^2 = \alpha x + \beta$, from which $1/x = (x - \alpha)/\beta$ is obtained.

A similar analysis is carried out for all other polynomial equations.

\subsection{Relationship with the arithmetical trapezoid}

The generating function of a particular sequence creates a power series whose coefficients are exactly the numbers that embrace such a sequence, following the same order built by the recurrence equation of the same sequence. For example, the power series $z + z^2 + 2z^3 + 3z^4 + 5z^5 + 8z^6 + 13z^7 + \cdots$, contains the numbers of the Fibonacci sequence, whose generating function is $f(z) = z/(1 - z - z^2)$. Many other generating functions and their respective power series can be found, with only a general function. In the case of a quadratic polynomial, the generating function is given as it follows
\begin{equation}
f(z) = \frac{x_0 + (x_1 - \alpha x_0)z}{1 - \alpha z - \beta z^2},
\end{equation}
where $\alpha$ and $\beta$  are , as previously, the polynomial coefficients and ($x_0$, $x_1$), also as previously, the sequence seeds. Thus, for the case of a cubic polynomial, there will be,
\begin{equation}
f(z) = \frac{x_0 + (x_1 - \alpha x_0)z + (x_2 - \alpha x_1 - \beta x_0)z^2}{1 - \alpha z - \beta z^2 - \gamma z^3},
\end{equation}
where, for this case, $x_0$, $x_1$ and $x_2$, are the sequence seeds, and $\alpha$, $\beta$ and $\gamma$, the polynomial coefficients.

In case the polynomial is of any $n$-degree, the generating function is given as,
\begin{widetext}
\begin{equation}
f(z) = \frac{\displaystyle\sum_{i = 1}^{n}\left[x_{n - i} - \Omega_{n - i - 1}\cdot\sum_{j = 0}^{n - i - 1}\alpha_{n - 1 - j}\,x_{n - i - 1 - j}\right]z^{n - i}}{1 - \displaystyle\sum_{i = 1}^{n}\alpha_{n - i}z^i}
\label{geradorageral}
\end{equation}
\end{widetext}
where $\Omega_n$ represents the unit function, defined as
\begin{equation}
\Omega_n = \left\{
\begin{array}{cc}
1 & \ ; \ (n \geq 0)\\[0.3cm]
0 & \ ; \ (n < 0)
\end{array}\right.
\end{equation}
and, $x_i$, the sequence seeds from the $n$-degree polynomial, and the $\alpha$'s, are the polynomial coefficients. This unit function is different from the conventional unit one. In the present case, the function will always be unitary for either any positive value or zero from $n$. In addition, it will be null for $n$ negative values.

It is with the generating function that the \textit{arithmetic trapezoids} can be built, and the steps for its construction will be shown below. Just for simplifying,   the numerator of Eq. (\ref{geradorageral}) will be represented as $T \equiv T(z)$ and the sum of the denominator as $R \equiv R(z)$. Thus, the equation will have the simple form $f(z) = T(z)/[1 - R(z)]$. By expanding the denominator of the expression in a Taylor's series, in terms of $R$, a product series  $T\cdot R^i$, will be obtained, each producing polynomials, that is,
\begin{equation}
\begin{array}{ccl}
T\cdot R^0 & = & x_0 + (x_1 - \alpha x_0) z + \cdots
\end{array}\nonumber
\end{equation}
\begin{equation}
\begin{array}{ccl}
T\cdot R^1 & = & \alpha x_0 z + [\alpha x_1 + (\beta - \alpha^2) x_0 + \cdots] z^2\\[0.1cm]
                  &     & \ \ + (\beta x_1 - \alpha\beta x_0 + \cdots) z^3 + \cdots
\end{array}\nonumber
\end{equation}
\begin{equation}
\begin{array}{ccl}
T\cdot R^2 & = & \alpha^2x_0z^2 + [\alpha^2x_1+(2\alpha\beta - \alpha^3)x_0 + \cdots]z^3\\[0.1cm]
                  &     & \ \ +[2\alpha\beta x_1+(\beta^2 - 2\alpha^2\beta)x_0 + \cdots]z^4\\[0.1cm]
                  &     & \ \ + (\beta^2x_1 - \alpha\beta^2x_0 + \cdots)z^5 + \cdots
\end{array}\nonumber
\end{equation}
\begin{equation}
\begin{array}{ccl}
T\cdot R^3 & = & \alpha^3x_0z^3 +[\alpha^3x_1+(3\alpha^2\beta - \alpha^4)x_0 + \cdots]z^4\\[0.1cm]
                  &     & \ \ + [3\alpha^2\beta x_1+(3\alpha\beta^2 - 3\alpha^3\beta)x_0 + \cdots]z^5\\[0.1cm]
                  &     & \ \ + [3\alpha\beta^2x_1+(\beta^3 - 3\alpha^2\beta^2)x_0 + \cdots]z^6\\[0.1cm]
                  &     & \ \ + (\beta^3x_1 - \alpha\beta^3x_0 + \cdots)z^7 + \cdots\\
                  &  \vdots  &
\end{array}\nonumber
\end{equation}
The coefficients of these polynomials will build the trapezoid; and the power $i$ value indicates the trapezoid row number. In this study only the trapezoid of the quadratic and cubic equations will be shown, that is, for $n = 2$ and $n = 3$, in Eq. (\ref{geradorageral}). The general case for the polynomial of $n$-degree will not be exhibited because it is an extremely extensive expression, but it is indicated as a further exercise for the reader.

By denominating the coefficients as $C_{i,j}$, the trapezoid for these two cases is built,
\begin{widetext}
\begin{center}
\hspace{2.cm}Quadratic Polynomial Trapezoid\\[0.3cm]
\begin{tabular}{ccccccccccccccc}
$T\cdot R^0$ &\ \ \ row 0 \ \ & & & & & & $C_{0,0}$ & & $C_{0,1}$ & & & & & \\
$T\cdot R^1$ &\ \ \ row 1 \ \ & & & & & $C_{1,0}$ & & $C_{1,1}$ & & $C_{1,2}$ & & & & \\
$T\cdot R^2$ &\ \ \ row 2 \ \ & & & & $C_{2,0}$ & & $C_{2,1}$ & & $C_{2,2}$ & & $C_{2,3}$ & & & \\
$T\cdot R^3$ &\ \ \ row 3 \ \ & & & $C_{3,0}$  & & $C_{3,1}$ & & $C_{3,2}$ & & $C_{3,3}$ & & $C_{3,4}$ & & \\
$T\cdot R^4$ &\ \ \ row 4 \ \ & & $C_{4,0}$ & & $C_{4,1}$ & & $C_{4,2}$ & & $C_{4,3}$ & & $C_{4,4}$ & & $C_{4,5}$ & \\
$T\cdot R^5$ &\ \ \ row 5 \ \ &\ \ \ $C_{5,0}$ & & $C_{5,1}$  & & $C_{5,2}$ & & $C_{5,3}$ & & $C_{5,4}$ & & $C_{5,5}$ & & $C_{5,6}$
\end{tabular}
\end{center}${}$\\
\begin{center}
Cubic Polynomial Trapezoid\\[0.3cm]
\begin{tabular}{cccccccccccccccccccccccccc}
row 0 \ \ &&  & & & & & & & & & $C_{0,0}$ & & $C_{0,1}$ &  & $C_{0,2}$ &  &  &  &  &  &  &  &  &  &   \\
row 1 \ \ && & & & & & & & $C_{1,0}$ & & $C_{1,1}$ &  & $C_{1,2}$ &  & $C_{1,3}$ &  & $C_{1,4}$ &  &  &  &  &  &  &  &   \\
row 2 \ \ && & & & & & $C_{2,0}$ & & $C_{2,1}$ &  & $C_{2,2}$ &  & $C_{2,3}$ &  & $C_{2,4}$ &  & $C_{2,5}$ &  & $C_{2,6}$ &  &  &  &  &  &   \\
row 3 \ \ && & & & $C_{3,0}$ & & $C_{3,1}$ &  & $C_{3,2}$ &  & $C_{3,3}$ &  & $C_{3,4}$ &  & $C_{3,5}$ &  & $C_{3,6}$ &  & $C_{3,7}$ &  & $C_{3,8}$ &  &  &  &   \\
row 4 \ \ && & $C_{4,0}$ & & $C_{4,1}$ &  & $C_{4,2}$ &  & $C_{4,3}$ &  & $C_{4,4}$ &  & $C_{4,5}$ &  & $C_{4,6}$ &  & $C_{4,7}$ &  & $C_{4,8}$ &  & $C_{4,9}$ &  & $C_{4,10}$ &  &   \\
row 5 \ \ &$C_{5,0}$ & & $C_{5,1}$ &  & $C_{5,2}$ &  & $C_{5,3}$ &  & $C_{5,4}$ &  & $C_{5,5}$ &  & $C_{5,6}$ &  & $C_{5,7}$ &  & $C_{5,8}$ &  & $C_{5,9}$ &  & $C_{5,10}$ &  & $C_{5,11}$ &  & $C_{5,12}$  
\end{tabular}
\end{center}
\end{widetext}
If the quartic equation is considered, the first row of its trapezoid, $i = 0$, will have four components, that is, $\{C_{00},C_{01},C_{02},C_{03}\}$, and the second row, $i = 1$, seven components,, $\{C_{10},C_{11},C_{12},C_{13},C_{14},C_{15},C_{16}\}$, and so on. In general, the trapezoid that corresponds to any polynomial of $n$-degree, will have in row $i$, $p = i\cdot(n - 1) + n$ components. Therefore, as the degree of the polynomial increases, the trapezoid is extensive.

The components $C_{ij}$ for these two cases, that is, the quadratic and cubic ones, are calculated with the following expressions,
\begin{widetext}
\begin{equation}
C_{ij} = \displaystyle\alpha^{i - j}\beta^{j - 1}\left\{\left[\binom{i}{j}\beta - \binom{i}{j - 1}\alpha^2\right]x_0 + \binom{i}{j - 1}\alpha x_1\right\},
\label{trapeziosegundograu}
\end{equation}
\end{widetext}
considering the quadratic polynomial case.  The factor $\binom{\ast}{\ast}$ represents Newton's binomial. In relation to the cubic polynomial case, there will be,
\begin{widetext}
\begin{equation}
\begin{array}{ccl}
C_{ij} &=& \displaystyle\alpha^{i - j}\sum_{k = 0}^{m}
\alpha^k\beta^{j - 2 - 2k}\gamma^k\left\{
\binom{j - k - 2}{k}\binom{i}{j - k - 2}\,\alpha^2\,x_2\right.\\[0.6cm]
&&\ \ \ \ \ \displaystyle +\ \alpha\left[\binom{j - k - 1}{k}\binom{i}{j - k - 1}\,\beta
- \binom{j - k - 2}{k}\binom{i}{j - k - 2}\,\alpha^2\right]x_1\\[0.6cm]
&&\ \ \ \ \ \ \displaystyle\left.+\ \beta\left[\binom{j - k}{k}\binom{i}{j - k}\,\beta
- \binom{i - k + 1}{j - 2k - 1}\binom{i}{i - k}\,\alpha^2\right]x_0
\right\}
\label{trapezioterceirograu}
\end{array}
\end{equation}
\end{widetext}
where
\begin{equation}
m = \delta_{0,j\%2}\cdot\left(\frac{j}{2}\right) + \delta_{0,(j + 1)\%2}\cdot\left(\frac{j - 1}{2}\right)
\end{equation}
where $\alpha$ and $\beta$ are the quadratic polynomial coefficients, and $\alpha$, $\beta$ and $\gamma$, are the cubic polynomial coefficients, already shown in Eqs. (\ref{segundograu}) and (\ref{terceirograu}), respectively. In addition $\delta_{ij}$ represents Kronecker's delta. The Fibonacci trapezoid $(x_0,x_1) = (0,1)$, which corresponds to the Fibonacci sequence, is as it follows,
\begin{center}
\begin{tabular}{ccccccccccccc}
& & & & & 0 & & 1 & & & & & \\
& & & & 0 & & 1 & & 1 & & & & \\
& & & 0 & & 1 & & 2 & & 1 & & & \\
& & 0  & & 1 & & 3 & & 3 & & 1 & & \\
& 0 & & 1 & & 4 & & 6 & & 4 & & 1 & \\
0 & & 1 & & 5 & & 10 & & 10 & & 5 & & 1
\end{tabular}
\end{center}
This trapezoid corresponds exactly to Pascal's triangle, if the `zero layer' is eliminated.

The trapezoid for another particular Fibonacci set $(x_0,x_1) = (2,1)$, which corresponds to Lucas' sequence, is as it follows,
\begin{center}
\begin{tabular}{ccccccccccccc}
& & & & & 2 & & -1 & & & & & \\
& & & & 2 & & 1 & & -1 & & & & \\
& & & 2 & & 3 & & 0 & & -1 & & & \\
& & 2  & & 5 & & 3 & & -1 & & -1 & & \\
& 2 & & 7 & & 8 & & 2 & & -2 & & -1 & \\
2 & & 9 & & 15 & & 10 & & 0 & & -3 & & -1
\end{tabular}
\end{center}
which can be called `Lucas' triangle' in case the trapezoid 'second layer'  is eliminated, although Lucas has never created this triangle.

In the silver case $(\alpha,\beta) = (2,1)$ and $(x_0,x_1) = (0,1)$, which corresponds to Pell's sequence, the trapezoid is as it follows\\
\begin{center}
\begin{tabular}{ccccccccccccc}
& & & & & 0 & & 1 & & & & & \\
& & & & 0 & & 2 & & 1 & & & & \\
& & & 0 & & 4 & & 4 & & 1 & & & \\
& & 0 & & 8 & & 12 & & 6 & & 1 & & \\
& 0 & & 16 & & 32 & & 24 & & 8 & & 1 & \\
0 & & 32 & & 80 & & 80 & & 40 & & 10 & & 1
\end{tabular}
\end{center}

Considering the cubic cases, the Fibonacci set $(x_0,x_1,z_2) = (0,1,1)$, known as tribonacci sequence, whose numbers are,
\[0, 1, 1, 2, 4, 7, 13, 24, 44, 81, 149, \dots\]
corresponds to the following trapezoid
\begin{center}
\begin{tabular}{ccccccccccccccccccccccccc}
   &&    &&    &&      &&      &&   0 &&   1 &&   0 &&      &&      &&    &&    &&   \\
   &&    &&    &&      &&   0 &&   1 &&   1 &&   1 &&   0 &&      &&    &&    &&   \\
   &&    &&    &&   0 &&   1 &&   2 &&   3 &&   2 &&   1 &&   0 &&    &&    &&   \\
   &&    && 0 &&   1 &&   3 &&   6 &&   7 &&   6 &&   3 &&   1 && 0 &&    &&   \\
   && 0 && 1 &&   4 && 10 && 16 && 19 && 16 && 10 &&   4 && 1 && 0 &&   \\
0 && 1 && 5 && 15 && 30 && 45 && 51 && 45 && 30 && 15 && 5 && 1 && 0  
\end{tabular}
\end{center}
And so on. Therefore, infinite trapezoids can be built by varying the values of both, the polynomial and seed coefficients.

Another important issue is that the outer layer of the left hand side of the trapezoid will only repeat the seed value $x_0$, regardless the values that the polynomial coefficients might have. When the outer layer is zero, it can be smoothly eliminated. In the quadratic case, the shape of a triangle will appear; it is the case of Pascal's triangle, when $(x_0,x_1) = (0,1)$. Considering the tribonacci case, with the elimination of the `zero layer' when $x_0 = 0$, the trapezoid will sometimes result in a triangle, but sometimes not. Thus, in general, it is better to be called trapezoid, including the corresponding Pascal.

Commonly, for an $n$-degree polynomial, a lot of interesting properties arise from the general trapezoid, generalizing the well-known Pascal's triangle ones. Among them, only three will be highlighted. One of them is analogous to the recurrence equation of expression (\ref{graugeral}),
\begin{equation}
C_{i+1,j+n-1} = \displaystyle\sum_{k = 0}^{n-1}\alpha_k C_{i,j+k}.
\end{equation}
where $n$ is the degree of the polynomial. This expression, when applied to a quadratic order Fibonacci, produces $C_{i+1,j+1} = C_{i,j+1} + C_{i,j}$, which is the usual rule of Pascal's triangle. The second is the sum of the components of a horizontal row,
\begin{widetext}
\begin{equation}
S_{i,n} = \sum_{j = 0}^{(i + 1)\cdot(n - 1)}C_{i,j} = \left[\displaystyle\sum_{r = 0}^{n - 1}\left(1 - \sum_{l = r}^{n - 2}\alpha_{n+r-l-1}\right)x_r\right]\cdot\left(\sum_{r = 0}^{n - 1}\alpha_r\right)^i
\end{equation}
\end{widetext}
where $i$ represents the trapezoid row number, and $n$, the degree of the polynomial. Once more, for a quadratic order Fibonacci, with $(x_0,x_1) = (0,1)$, the sum will be, $S_{i,2} = 2^i$. Finally, the implicit sequence number in the trapezoid was rescued by just adding the diagonal components of each row $i$, that is,
\begin{equation}
x_i = \sum_{j = 0}^i\, C_{i-j,j}.
\end{equation}
It should be highlighted that the relationship of this trapezoid with the numeric sequences is natural and the link between them is the generating function.

To conclude this topic, it can be said that a multitude of arithmetical trapezoids can be obtained, only with Eqs. (\ref{trapeziosegundograu}) and (\ref{trapezioterceirograu}), which are associated with the quadratic polynomial, (\ref{segundograu}), and the cubic one, (\ref{terceirograu}), and related to the numeric sequences (\ref{xn2grau}) and (\ref{xn3grau}).

When choosing a polynomial function of any $n$-degree by using Eqs. (\ref{graugeral}), (\ref{equacaolinear}) and (\ref{sistemageral}), the number of numeric sequences is infinite, and only a single arithmetic trapezoid corresponds to one of such sequences.\\

\section{Conclusion}

The aim of this study was to show the relationship among polynomials, numerical sequences and numerical trapezoids being the Pascal's Triangle a very particular case. The starting point, as seen before, is the polynomial equations. The coefficients of these polynomials generalize the numerical sequences and the golden numbers to which these sequences are linked. Comprehensibly, these golden numbers are nothing more than roots of the polynomial equations and attract such sequences. It was also shown that the arithmetic trapezoids are related to the numeric sequences, which in turn, are related to the polynomials; and that their construction is carried out through the generating functions.

Regarding trapezoids, one of the applications is inserted in the probability and distribution field. This is a wide area since such distributions are associated with non-extensiveness contexts. Considering this subject, further studies will be carried out.

\end{document}